\documentclass{amsart}[12pt]
\usepackage{bbm}
\usepackage{amsmath}
\usepackage{amssymb}
\usepackage{amsthm}
\usepackage{setspace}
\usepackage{amsfonts}
\usepackage{epsfig}

\theoremstyle{plain}
\newtheorem{theorem}{Theorem}[section]
\newtheorem{lemma}[theorem]{Lemma}
\newtheorem{claim}[theorem]{Claim}

\theoremstyle{definition}
\newtheorem{definition}[theorem]{Definition}

\title{Hyperbolic outer billiards : a first example}
\author{Daniel Genin}

\begin{document}

\begin{abstract}
We present the first example of a hyperbolic outer billiard.  More precisely
we construct a one parameter family of examples which in some sense
correspond to the Bunimovich billiards.
\end{abstract}
\maketitle

\section{Introduction}

Since the introduction of outer billiards by B.H. Neumann in 1959 \cite{Neu1}
and their popularization by J. Moser in \cite{Mo1} and \cite{Mo2}. there have
been several developments indicating that their dynamics in some respects
parallels  the dynamics of ordinary billiards \cite{Ta1}, \cite{Ta2},
\cite{Gu-Ka}, \cite{Bo1}.  For example, it has been shown that the analogue of
the Lazutkin's theorem holds for outer billiards \cite{Do1}, and that the string
construction in ordinary billiards has its parallel in the secant area
construction of outer billiards\cite{Be}.  In contrast, the rich theory of
chaotic billiards so far has no counterpart, in significant part due to the lack
of examples.  The aim of the present paper is to make a first step toward
eliminating this void, by providing a first example of a  hyperbolic outer
billiard.

More precisely we produce a one parameter family of outer billiards which have a
square invariant region of positive measure on which the Lyapunov exponents of
the outer billiard map are non-zero almost everywhere.  If one were to draw an
analogy with ordinary billiards this family of examples most closely parallels
the billiards of Bunimovich \cite{Bu1}, in that the billiards are composed of
arcs of hyperbolas which individually have "integrable" dynamics joined by
"neutral" segments -- corners.  As in the case of Bunimovich billiards
non-vanishing of the Lyapunov exponents guarantees positive entropy and other
nice properties.  Also as in the case of Bunimovich billiards ergodicity doesn't
come for free and has to be proved separately.  We expect, nevertheless, that
the outer billiards in our family of examples are indeed ergodic in the invariant
region.

Dynamics outside the invariant region appears to be non-hyperbolic.  Numerical
explorations indicate the presence of invariant curves and elliptic islands.  So
this family of outer billiards, in addition, provides a nice example of coexistence of
hyperbolic and elliptic behavior.

The proof of the main result uses the cone field method introduced by Wojtkowski
in \cite{Wo1}.  After producing an invariant region we construct a measurable
field of cones defined almost everywhere in the invariant region and then show
that the cones are eventually strictly preserved.  We begin the exposition with
a brief introduction to outer billiards followed by a section on the secant area
construction which allows us to obtain an invariant region, or to be exact a
table for a given invariant region; in the following section we define the cone
field over the invariant region, and then prove that it is eventually strictly
preserved in the last section.

\section{Outer billiards}
We begin with a definition of the outer billiard map

\begin{figure}[h]
\centering
\epsfig{file=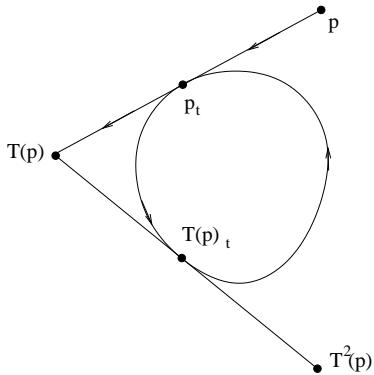,width=5cm}
\caption{Definition of outer billiard map}
\end{figure}

\begin{definition}
Let $\Gamma$ be an oriented strictly convex plane curve in $\mathbf{R}^2$ and $\Omega$
the domain enclosed by $\Gamma$.  Suppose for the moment that $\Omega$ is
strictly convex and let $D=\mathbf{R}^2\setminus \Omega$.  The \emph{outer
billiard map} is the continuous transformation $T:D\rightarrow D$ uniquely
specified by
\begin{itemize}
\item[(i)] the oriented segment $[p,T(p)]$ is tangent to $\Gamma$ at some $p_t$;
\item[(ii)] orientation of $\Gamma$ agrees with orientation of the segment $[p,T(p)]$ at $p_t$;
\item[(iii)] $p_t$ bisects $[p,T(p)]$
\end{itemize}
\end{definition}

In analogy with inner billiards we will say that $p$ reflects at the boundary in
$p_t$. It is easy to see that the resulting map is a continuous transformation
of $D$. The condition of strict convexity can be relaxed to allow outer
billiards with discontinuities.  Flat segments of the boundary are equivalent to
corners for inner billiards -- orbit of a point reflecting in a flat segment
can not be extended.  Outer billiards about polygons are of this type and have
been studied extensively in recent times as part of a more general program to
understand complexity arising in piece-wise isometries (\cite{K-L-V},
\cite{Gu-T}, \cite{Ka1}).

It is also easy to show that $T$ is an area preserving twist map, the invariant
measure being simply the Lebesgue measure on $D$ (\cite{Ta1}).  

Outer billiards have many remarkable properties (see for example \cite{Ta1},
\cite{Ta2}, \cite{Ko1}, \cite{S-V}, \cite{Gu-S},).

\section{Secant area construction}
We begin by describing a construction that given a convex plane
curve produces an outer billiard table for which the curve is an
invariant one.  This construction parallels the string
construction of inner billiards which does the same for caustics.

Given a convex plane curve $\gamma$ and a parameter $a$ satisfying $0<a<A$,
where $A$ is the area enclosed by $\gamma$, we consider the family of lines
$\mathcal{L}$ that divide the region enclose by $\gamma$ into parts with area $a$
and $A-a$. We have the following result

%
\begin{lemma}
The envelope of $\mathcal{L}$ is a closed curve $\Gamma$, which is convex if
there are no cusps. Furthermore, $\gamma$ is an invariant curve of outer
billiard about $\Gamma$.
\end{lemma}
%
The proof of this result can be found in \cite{Ta2}.

Using the area construction we can construct a table with an arbitrary convex
invariant region.
%
\section{Table construction}
\label{construction}
We consider a 1-parameter family of outer billiards obtained by
the area construction from a the unit square $\gamma$.
\begin{figure}[h]
\centering
\epsfig{file=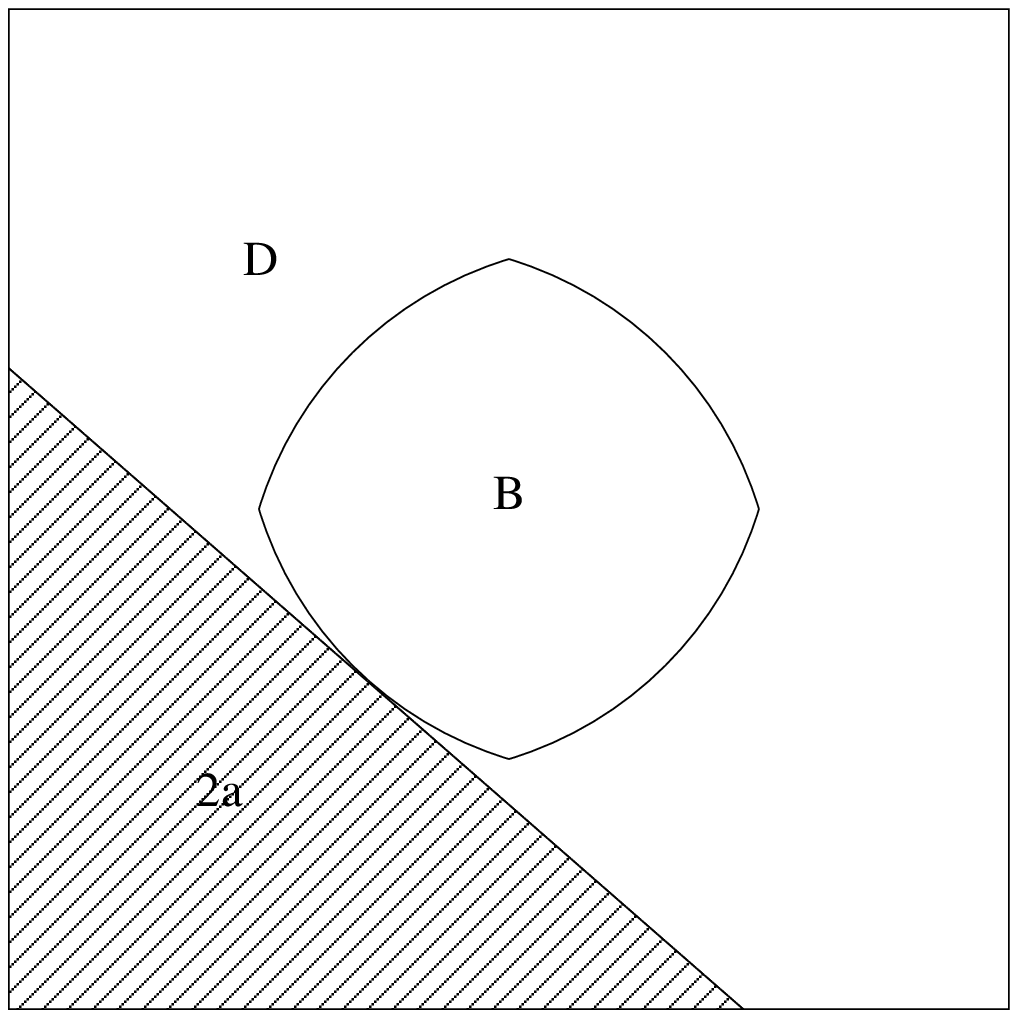,width=5cm}
\caption{The outer billiard table}
\end{figure}
Fix the area parameter to be $2a$ and let $B$ be the (closed) billiard
table, i.e. the closed region enclosed by $\Gamma$, and $D$ be the complement of $B$ in the open unit square. Let $\Gamma$ be oriented in the counter-clockwise direction.
%
\begin{lemma}
$\Gamma$ is composed of four arcs of hyperbolas, asymptotic to
lines containing the sides of $\gamma$, meeting in corners
opposite centers of sides of $\gamma$.
\end{lemma}
%
The proof of the above lemma is elementary and is left to the reader.
Depending on which arc of $B$ contains $p_t$, it is convenient to
use a different coordinate system on $D$.  Each arc naturally corresponds
to a corner, which we choose to be the origin when considering points
reflecting in this arc.  Choose the outgoing edge (considering $\gamma$
oriented clockwise) to be the $x$-axis and the incoming -- the $y$-axis. 
We introduce corresponding coordinate functions $x(p)$ and $y(p)$.  In these
coordinates the hyperbolic arcs of $B$ from the above lemma are given by the
equation $y(p)x(p)=a$.

As a consequence of the construction $D$ is an invariant region
for the outer billiard about $\Gamma$.
%
\section{Cone field}
We will define an invariant cone field on a subset $D_h$ of $D$
defined by 
$D_h=\{p\in D| 
T^k(p)_t \mathrm{is not a corner of}\ B \mathrm{for some}\ k\in\mathbf{Z}\}$.  
That is $D_h$ is the set of all points in $D$ whose orbits touch the interior of 
one of the hyperbolic arcs of $B$ at least once.
%
\begin{lemma}
\label{periodic}
The points in $D\setminus D_h$ are periodic.
\end{lemma}
\begin{proof}
A point in $D\setminus D_h$ reflects only in corners of $B$.  Hence as far as
it is concerned the table is a square.  But it is well known that outer
billiard about any lattice polygon has only periodic points (see for example
\cite{Ta1}).

\end{proof}
%

The set $D\setminus D_h$ is indicated in Figure \ref{Dh_diagram} in black.  It
is easy to see that all points in $D\setminus D_h$ have the same period --
four.  So $D\setminus D_h$ is fixed by $T^4$ and hence is an elliptic domain.

A cone $C(p)\subset T_pD_h$ will be defined by a pair of vectors $(u,v)\in
{T_pD_h}^2$, $C(p)=\{w\in T_pD_h|[u,w][w,v]>0\}$ where $[\cdot,\cdot]$ stands
for the cross product induced by identifying $T_pD_h$ with the ambient
$\mathbf{R}^2$.  It is clear that scaling $u$ and $v$ in the above definition by
positive constants does not change $C(p)$ so whenever we will talk about
equality of vectors defining cones we shall mean equality up to multiplication by
a positive constant.

We will first define the invariant cone field on the points in $D_h$ for which
$p_t$ is not a corner of $B$.  The cone field can then be extended to all
points in $D_h$ which do not land on discontinuities of $dT$ by pulling back the
cones by $dT$.  We set $C(p)=(u(p),v(p))$, where for $p$ with $x(p)=x$, $y(p)=y$
$$u(p)=(1,-y/x)$$
is a vector tangent to the homothetic hyperbola passing through
$p$, and
$$v(p)=p-p_t=\left(x-\frac{a+\sqrt{a^2-axy}}{y},y-\frac{a y}{a+\sqrt{a^2-axy}}\right)$$
In the last formula we used
$$(x(p_t),y(p_t))=\left(\frac{a+\sqrt{a^2-axy}}{y},\frac{a y}{a+\sqrt{a^2-axy}}\right)$$

\begin{figure}[h]
\centering
\epsfig{file=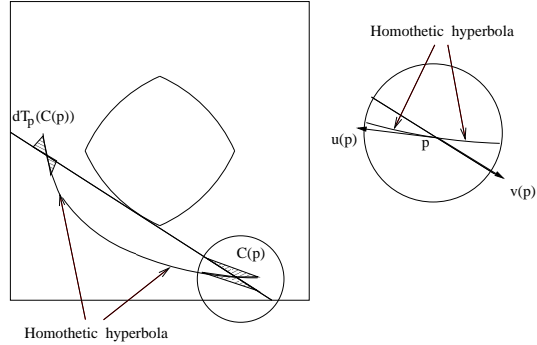,width=7cm}
\caption{Cone construction}
\label{cones}
\end{figure}

For convenience we will define the cone field at points of $\gamma$ (which are
not in $D$) by taking $u(p)$ to be the tangent vector to $\gamma$ in the
corresponding direction.
\section{Hyperbolicity}
Now we are ready to prove the main result.
%
\begin{theorem}
\label{hyperbolic}
$C(p)$ is eventually strictly preserved under $T$ a.e. in $D_h$ for 
$a\in (0,(3-\nobreak\sqrt{5})/8)$.  Hence $T$ has non-vanishing Lyapunov
exponents a.e. in $D_h$ for the corresponding set of parameter values.
\end{theorem}
%
We will prove this by showing that the cone field defined in the previous
section is eventually strictly preserved.  
\begin{definition}
\hspace{0ex}
\begin{enumerate}
\item A cone field $C(p)$ is \emph{preserved} at a point $p$ if $dT_p
C(p)\subset C(T(p))$.  

\item A cone field is \emph{strictly} preserved if 
$dT C(p)\subset \stackrel{\circ}{C}(T(p))$.

\item A cone field is \emph{eventually} strictly preserved if for almost every
$p$ there exists an $n(p)$ such that the cone field is strictly preserved at
$T^{n(p)}(p)$.
\end{enumerate}
\end{definition}
The result above then follows from the result of Wojtkowski introduced in
\cite{Wo2}.
%
\begin{theorem}
If there exists a measurable bundle of sectors which is eventually strictly
preserved by $\Phi:M^2\rightarrow M^2$, where $M^2$ is a two dimensional Riemannian
manifold, satisfying
\begin{enumerate}
\item $\Phi$ preserves a probabilistic measure $\mu$ which has a non-vanishing
density with respect to the Riemann area element on $M^2$;
\item The singularities of $\Phi$ satisfy
$$\int_{M^2} \log^+\|D_x \Phi^{\pm 1}\|d\mu(x)<+\infty$$
where $\log^+t=\max(\log t,0)$.
\end{enumerate}
Then the Lyapunov exponent $\lambda_+$ of $\Phi$ is positive $\mu$ a.e.
\end{theorem}
%
$T$ satisfies the first condition above because it preserves the Lebesgue
measure on $D_h$.  It also satisfies the second condition because the
differential of $T$ is bounded in norm on $D_h$.  Indeed, differential of an
outer billiard map blows up only if the table has points of vanishing
curvature, which $B$ does not.  Thus the above result may be applied to
$T$ if $C(p)$ is eventually strictly preserved.
This will be proved in a series of lemmas because the argument
naturally divides into cases.  The argument is essentially different
for points with $p_t$ and $T(p)_t$ belonging to the interior of the same
hyperbolic segment of $B$ and for points for which they belong to interiors of
different hyperbolic segments.

We begin by examining the first case.
\begin{lemma}
\label{same_segment}
$C(p)$ is preserved at points such that $p_t$ and $T(p)_t$ belong to the
interior of the same hyperbolic segment.
\end{lemma}
\begin{proof}
Since we are concerned with one hyperbolic arc we can for the
moment forget about the rest of the table and consider outer
billiard about a single branch of hyperbola.

Simple algebra shows that $xy$ is an integral of motion for $T$ and
hence the homothetic hyperbolas $xy=c$ are preserved by $T$.
Since $T$ preserves the order of points on the invariant hyperbolas the vector $u(p)$
tangent to a hyperbola at $p$, is mapped by $dT$ to a vector tangent to the same hyperbola at $T(p)$.  That is $dT_p u(p)=u(T(p))$. Also the tangent vector 
$v(p)=p-p_t$ is mapped to minus itself $dT_p v(p)=-v(p)$ since restricted to the line of
tangency the billiard map is simply a reflection in the tangency point.
Hence the image cone $dT_p C(p)=(u(T(p)),-v(p))$.  Noting that the vectors
$u(T(p))$ and $dT_p u(p)$ are the same up to rescaling by a positive constant it
is enough to check that $dT_p v(p)\in C(T(p))$.  We will show that the angle
of $C(p)$ is always obtuse and the angle of $dT_p C(p)$ is always acute. 
Observe that $v(p)$ always lies in the fourth quadrant while $u(p)$ always lies
in the second.  So the angle between them is always obtuse.  
Similarly $u(T(p))$ is always in the second quadrant and so is
$-v(p)$ so the angle between them is always acute.  Hence 
$dT_p C(p)\subset C(T(p))$ for points reflecting in a single hyperbolic segment.

Note that $dT_p v(p)$ is strictly inside $C(T(p))$ although
$dT_p u(p)=\nobreak u(T(p))$ up to rescaling by positive constant so the inclusion is not
strict.
\end{proof}
%
The following consequence of the proof is useful in itself.
%
\begin{lemma}
\label{angles}
For every point $p\in D_h$ such that $p_t$ is in the interior of a
hyperbolic arc.
\begin{enumerate}
\item $C(p)$ has an obtuse angle
\item $dT_p C(p)$ has an acute angle
\end{enumerate}
\end{lemma}
%
Before moving on to the second case we observe that cones $C(p)$ are
nested in a special way along the line through $p$, $p_t$.

\begin{figure}[h]
\centering
\epsfig{file=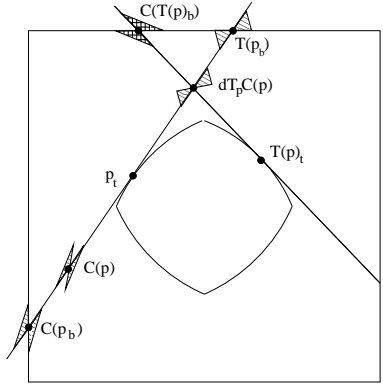,width=5cm}
\caption{Cone nesting}
\label{nesting}
\end{figure}

\begin{lemma}
Let $p$ and $p'$ be two points satisfying $p_t=p'_t$ and
$|p-\nobreak p_t|>|p'-p_t|$ then
\begin{enumerate}
\item $C(p)\subset C(p')$
\item $dT_{p'} C(p')\subset dT_p C(p)$
\end{enumerate}
where the tangent spaces are identified by parallel translation to
determine the inclusion relationship between the cones (Figure \ref{nesting}).
\end{lemma}
%
Proof of this statement involves only elementary geometry and is left to the 
reader.
The above lemma provides an easy way of determining whether a cone
field at $p$ is preserved.  Let $p_b$ be the intersection of the
ray from $p_t$ through $p$ with $\gamma$, then to check that the
cone field is preserved at $T(p)$ it is enough to check the inclusion
between cones at $p_b$ and $T(p)_b$. Indeed, if
$dT_{p_b}C(p_b)\subset C(T(p)_b)$ then
$$dT_pC(p)\subset dT_{p_b}C(p_b) \subset C(T(p)_b) \subset C(T(p))$$
so
$$dT_pC(p)\subset C(T(p))$$
Note that for $p\in \stackrel{\circ}{D}_h$ if $dT_{p_b}v(p_b)$ is in the interior of
$C(T(p)_b)$ then $dT_pC(p)\varsubsetneq\nobreak C(T(p))$, i.e. the inclusion
becomes strict, because from the above lemma $dT_p u(p)$ is in the interior of
$dT_{p_b} C(p_b)$ (Figure \ref{nesting}).

We now proceed to examine the points for which $p_t$ and $T(p)_t$ belong
to interiors of different hyperbolic segments of $B$ or corners.

%
\begin{lemma}
$C(p)$ is strictly preserved at points such that $p_t$ and $T(p)_t$ belong to
interiors of different hyperbolic segments.
\end{lemma}
\begin{proof}

If $p_t$ and $T(p)_t$ belong to interiors of different hyperbolic segments these
segments must be adjacent.  In this case (Figure \ref{nesting}) $T(p_b)$
and $T(p)_b$ belong to the same side of $\gamma$ which means that
$dT_{p_b}u(p_b)=u(T(p)_b)$.  Furthermore, $dT_{p_b}v(p_b)=v(T(p_b))\in
\stackrel{\circ}{C}(T(p)_b)$ by Lemma \ref{angles}. This implies that $dT_p(C(p))\varsubsetneq C(T(p))$. 

\end{proof}

\begin{lemma}

$C(p)$ is strictly preserved at points such that $p_t$ is a corner point.

\end{lemma}
\begin{proof}
Considering all points such that $p_t$ is a corner point we
obtain Figure \ref{Dh_diagram} which shows points reflecting in corners and
their images.

\begin{figure}[h]
\centering
\epsfig{file=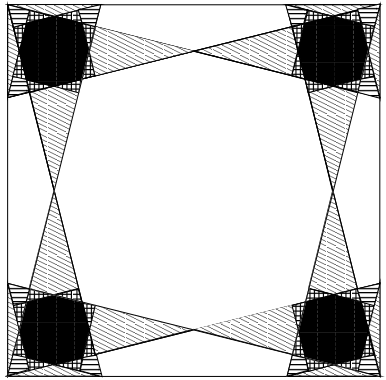,width=5cm}
\caption{Orbits reflecting in corners}
\label{Dh_diagram}
\end{figure}

In Figure \ref{Dh_diagram} all points reflecting in corners are shaded : those
that reflect in one corner between reflections in the interiors of sides are
shaded diagonally, those that reflect in two corners are shaded horizontally,
three -- doubly shaded.  Any point that hits all four corners is periodic, as
was proved in Lemma \ref{periodic}, and these points are solidly shaded in
black.  The lines in the diagram are the lines of discontinuity of the
derivative and their preimages.  We will use this diagram to reduce verification
of preservation of the cone field at shaded points to consideration of orbits of
a few representative points.  For this it is important to note that the diagram
in Figure \ref{Dh_diagram} is structurally invariant for  $a\in (0,1/4)$ (which
includes the parameter range under consideration), the interval in which the
lines tangent to $B$ and passing through the corners of $\gamma$ do not
intersect in the interior of $D$.

We first argue that we only need to consider preservation inside a few of the
shaded polygons.  This follows because the cone field and the map are both
preserved under the action of the rotation subgroup of $D_4$.  Hence if the cone
field is preserved inside some shaded polygon of Figure \ref{Dh_diagram} then
it is preserved in every polygon of its orbit under the rotation subgroup.
Therefore, it is enough to check preservation for one polygon per orbit.

Next, we show that for points in the shaded polygons only cones at the first
and last points of the orbit segment inside the shaded set need to be
compared.  $T$ is a central symmetry for every point $p$ with $p_t$ a corner
of $B$.  So $dT_p$ at such a point is $-I$ and $dT_pC(p)=C(p)$.  Hence it is
enough to check inclusion of cones at successive points of the orbit that
reflect in the interiors of sides.  Thus we ignore the points of the
orbit segment that reflect in corners, and let $p$ and $p'$ be the first and
last points respectively such that $p_t$ and $p'_t$ belong to the interiors of
sides.  As before $p_b$ and $p'_b$ will be intersections of the
corresponding rays with $\gamma$.

We further show that for points in the interior of a given polygon of the above
diagram there is only one way for the boundary points $p_b$, $T(p)_b$,...,
$p'_b$ to be distributed on the sides of $\gamma$.  This allows us to determine
preservation of cones for points in a given polygon by examining a diagram like
Figure \ref{nesting} for one of its interior points.  More precisely
Suppose the sides of $\gamma$ are numbered (exactly how is not important).  Then
to every point $p$ in one of the shaded polygons we can assign a sequence of
numbers of length at most 4 such that if the $k$-th symbol of the sequence is
$j$ then $T^k(p)_b$ belongs to the side of $\gamma$ labeled by $j$.
\begin{claim}
The sequence described above is the same for every point in a given shaded
polygon.
\end{claim}

\begin{proof}
Suppose there are two points $p$ and $q$ for which $T^k(p)_b$ and $T^k(q)_b$ lie
on different sides for some $k$.  Since the polygons are convex $T^k(p)$ and
$T^k(q)$ can be joined by a line segment contained in the interior of the
polygon.  By continuity there will be a point $r$ on this line segment such that
$r_b$ will be a corner point.  This is a contradiction since lines tangent to
$B$ and passing through corners of $\gamma$ form the boundaries of the shaded
polygons.
\end{proof}

In what follows we will say that a shaded polygon (or point) has order $k$ if it
belongs to an orbit segment that undergoes $k$ reflections in corners between
successive reflections in interiors of sides.  We will consider polygons of each
order in turn. 

\vspace{1ex}
\noindent I) We first consider order one domains.  Here there are two possibilities
corresponding to two different kinds of order one domains in the above diagram
(see Figure \ref{one_corner}).
\begin{figure}[h]
\centering
\epsfig{file=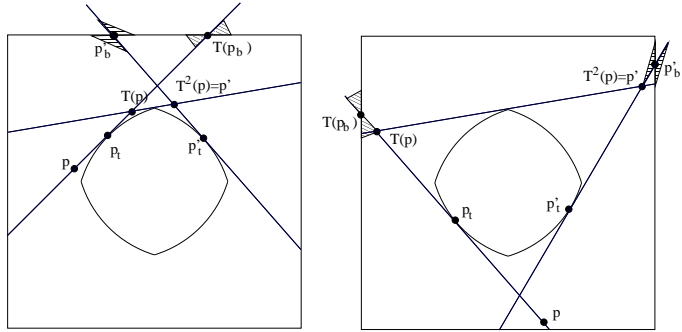,width=9cm}
\caption{Order one orbits}
\label{one_corner}
\end{figure}
In this case, it is enough to look at which sides the cones lie on to determine
inclusion. In the first case (the left diagram in Figure \ref{one_corner}), the
situation is identical to the case of a point touching adjacent
hyperbolic segments, and so the cone field is strictly preserved for these
points.  In the second case the cones $dT_{p_b}C(p_b)$ and $C(p'_b)$ lie on
opposite sides so $u(p'_b)=-dT_{p_b}u(p_b)$. As before, because the angle
between $dT_{p_b}u(p_b)$ and $dT_{p_b}v(p_b)$ is always acute and the angle
between $u(p'_b)$ and $v(p'_b)$ is always obtuse (Lemma \ref{angles}), 
$dT_{p_b}v(p_b)\in C(p'_b)$ and we have the desired strict inclusion.

\vspace{1ex}
\noindent II) For points touching two corners there is only one possible configuration
because the restriction on $a$ in the statement of
the theorem guarantees that a point can touch only neighboring corners of $B$. 
Indeed, if there was a point that reflected in opposite corners then one easily
checks that $a$ must be greater than $1/4$ which is ruled out by assumption.
Hence the only possible configuration is as in Figure \ref{two_corners}.
\begin{figure}[h]
\centering
\epsfig{file=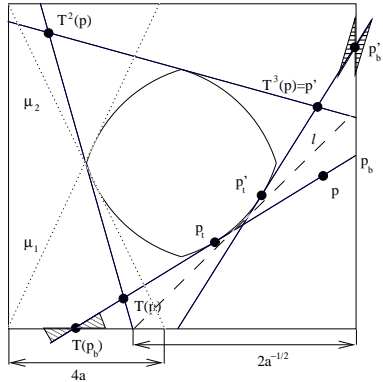,width=5cm}
\caption{Order two orbit}
\label{two_corners}
\end{figure}

The cone is strictly preserved if $p'_b$ precedes $p_b$ (the order is given by
the clockwise orientation of the square).  From the diagram it is clear that
this is true when $T(p)$ and $p'$ are on the same side of the line
$\ell$, tangent to the midpoint of the hyperbolic segment containing $p_t$ and
$p'_t$.  Since $T(p)$ reflects in a corner it has to lie between the
tangents to the hyperbolic arcs that meet at that corner, call them $\mu_1$ and
$\mu_2$ which are in turn lines in $\mathcal{L}$ that pass through corners of
$\gamma$.  The same reasoning applies to $p'$ for reverse time since $T^2(p)$
also reflects in a corner.  The above condition will be satisfied if
the triangle bounded by $\mu_1$, $\mu_2$ and $\gamma$, containing $T(p)$ does
not intersect $\ell$.  In this case $T(p)$ is guaranteed to lie on the right
side of $\ell$.  One easily checks that this happens exactly when $2\sqrt{a}+4a<1$
or $a<(3-\sqrt{5})/8$ as in the statement of the theorem.
\begin{figure}[h]
\centering
\epsfig{file=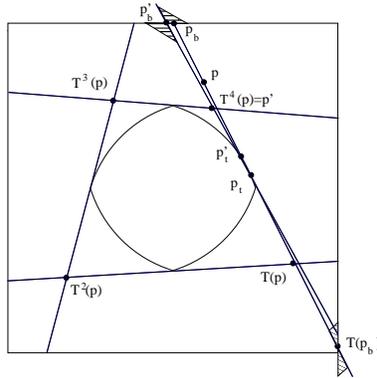,width=5cm}
\caption{Order three orbit}
\label{three_corners}
\end{figure}

\vspace{1ex}
\noindent III) Finally, for points touching three corners there is one possible
configuration as in Figure \ref{three_corners}.
From the diagram we see that the cone field is strictly preserved if
$p'_b$ precedes $p_b$.  Observe that reflecting a point in three consecutive
corners of $B$ gives a central symmetry about the fourth corner since the
corners of $B$ are vertices of a square.  It follows that $p'$ and $T(p)$ are
symmetric about the one corner untouched by the orbit segment from $p$ to $p'$,
and so that $[p',T(p)]$ contains this corner.  On the other hand $[p',T(p)]$
must intersect the interior of $B$ because otherwise $p'_t$ is a corner and
then $p$ has order four and so is not in $D_h$.  So the line through $p$ and
$T(p)$ must leave $p'$ on the same side as $B$.  
Noting that $p_t$ and $p'_t$ belong to the interior of
the same hyperbolic arc it follows that $p'_t$ precedes $p_t$, if the boundary
of $B$ is oriented clockwise, and so $p'_b$ precedes $p_b$.

\vspace{1ex}
We have proved that the cone field is strictly preserved at every point of
$D_h$ with the exception of points that land on discontinuities of the
derivative which form a set of measure zero.  Non-vanishing of the Lyapunov
exponents now follows.
\end{proof}

\section{Conclusion}

The above result shows that chaotic behavior is possible for outer billiards.
The numerical studies \cite{Ge1} indicate that it coexists with KAM-type
behavior on the rest of the domain and that near-integrable behavior may still
persist in the neighborhood of infinity.  It also appears that this one
parameter family of examples is a member of a much larger class of chaotic outer
billiards.  This class, we expect, contains many other outer billiards obtained
by the secant area construction from polygons.  Some numerical studies along
these lines for outer billiards obtained from a regular pentagon and a particluar
non-regular hexagon are also contained in the previous reference.  Unfortunately, 
the current proof does not appear to extend easily to these potential candidates.

\section{Acknowledgement}
I am grateful to my advisor Sergei Tabachnikov for setting me on this problem
and for many valuable discussions.

\clearpage
\bibliographystyle{plain}
\bibliography{article}

\end{document}